\documentclass[12pt]{report}
\usepackage{amssymb,amsmath,makeidx,verbatim}
\newcommand{\R}{\mathbb{R}}

\newcommand{\N}{\mathbb{N}}
\newcommand{\C}{\mathbb{C}}
\newcommand{\Z}{\mathbb{Z}}

\newcommand{\be}{\begin{enumerate}}
\newcommand{\ee}{\end{enumerate}}
\newcommand{\bq}{\begin{eqnarray*}}
\newcommand{\eq}{\end{eqnarray*}}

\begin{document}
%\pagenumbering{roman}
\newcommand{\disp}{\displaystyle}
\thispagestyle{empty}
\begin{center}
\textsc{Galois groups of Fermat polynomials and the arithmetic groups of Diophantine curves.\\}
\ \\
\textsc{Olufemi O. OYADARE}\\
\end{center}
\baselineskip 4pt
\begin{quote}
{\bf This paper develops a framework of algebra whereby every Diophantine equation is made quickly accessible by a study of the corresponding row entries in an array of numbers which we call the \textit{Newtonian triangles.} We then apply the framework to the discussion of some notable results in the theory of numbers. Among other results, we prove a new and complete generation of \textit{all} Pythagorean triples (without necessarily resorting to their production by examples), convert the collection of Newtonian triangles to a Noetherian ring (whose (multiplicative) identity element is found to be the Pascal triangle) and develop an easy understanding of the \textit{original} Fermat's Last Theorem ($FLT$). The application includes the computation of  the Galois groups of those polynomials coming from our outlook on $FLT$ and an approach to the explicit realization of arithmetic groups of curves by a treatment of some Diophantine curves.}
\end{quote}

\baselineskip 2pt
\ \\
\indent {\bf \S 1.  Introduction.} Let $x,y, n\in \mathbb{N} \cup\{0\},$ then the coefficients in the expansion of $(x+y)^{n},$ when considered as a polynomial in descending powers of $x,$ are $1,\;^nC_1y,\;^nC_2y^2, \cdots,y^n.$ For $y=1$ these coefficients form the $n$th row of the Pascal triangle, while, for other values of $y,$ the coefficients form the $n$th row of an array of numbers which we call the \textit{Newtonian triangles.} Numbers formed from these coefficients, by the application of the \textit{digital-correspondence} map, are $n-$powers of natural numbers and may be extended to generate all $n-$powers of rational numbers only. This outlook simplifies every Diophantine equation and gives proof of results that are consistent with the expectations of their originators and true to the spirit of classical number theory, as we shall show in the case of rational solutions of the equation $u^{n}+v^{n}=w^{n},$ for $n=2$ and its impossibility for non-zero rationals $u$, $v$ and $w,$ when integers $n>2,$ in $\S 3.$ and $\S4.,$ respectively.\\

\textit{The ideas of this paper emanated from a very elementary transformation of the finite Binomial theorem.} After the introduction of the digital-correspondence map and the Newtonian triangles in \S$2,$ we state and establish a purely algebraic reason for the existence and explicit form of all rational Pythagorean triples, leading to the partitioning of the integral ones in \S $3.$ Aside other mentioned approaches that may be taken to the study of Pythagorean triples, the ring of Newtonian triangles is introduced and proved to be Noetherian. \S $4$ contains an elementary proof of the \textit{original} Fermat's Last Theorem which is seen to be greatly simplified by the introduction and investigation of some built-in polynomials of the Newtonian triangles. Open problems on the ideal theory of the Noetherian ring of Newtonian triangles, distribution and density of solutions of Diophantine equations, non-rational Pythagorean triples in other \textit{fields} and the link with the \textit{Wiles-Taylor proof} of $FLT$ are all brought up in the remark at the end of each section. \S$5.$ contains two Lemmas and a Theorem, on the nature of those polynomials we call \textit{Fermat polynomials,} while we offer a novel approach to the yet-to-be-solved problem of computing the \textit{Mordell-Weil} groups of algebraic curves in \S6. Some open problems are contained in \S7.\\
\ \\

A preliminary version of Theorem $3.1$ is contained in the announcement $[9.].$\\

\ \\
\ \\
\ \\
$\overline{2010\; Mathematics}$ \textit{Subject Classification:} $11G15, \;\; 11Y40, \;\; 11Rxx, \;\; 14Kxx$
\ \\
\ \\
Uploaded on $19$ th May, $2014.$\\
\ \\
\ \\

\indent {\bf \S 2.  Digital-correspondence and Newtonian triangles.} A typical row in the Pascal triangle is $(1,\;^nC_1,\;^nC_2, \cdots,1).$ Among its properties we have that $1+^nC_1+^nC_2+ \cdots+1=2^{n},$ for all $n \in \N \cup\{0\}.$ For $n<5,$ each of the coefficients $1,\;^nC_1,\;^nC_2, \cdots,1,$ is a digit, so that any row may be viewed as a number having these coefficients as its digits. These numbers are $1,11,121,1331,$ and $14641,$ each of which is the respective $n$th power of $11,$ for $n=0,1,2,3,4.$ ($[10.],\;p.\;10$) It may then be asked:\\

\textit{Is it a mere coincidence that for $n \in \N \cup \{0\},\;n<5,$ the number $(1+1)^{n}$ (where $1$ is the repeated digit of the number $11=(10+1)$) is exactly $2^{n}$ (the sum $1+^nC_1+^nC_2+ \cdots+1$)? Indeed, what can we say of each of the remaining rows in the Pascal triangle with respect to $(11)^{n}$?}\\
\ \\

We answer the second question above as follows. Since the $5$th row in the triangle is $(1,5,10,10,5,1)$ an appropriate transfer of tens, at the middle terms, gives the number $161051.$ This is $11^{5}.$ We have taken the top digit $1$ in the Pascal triangle as the $0$th row. The $6$th row is $(1,6,15,20,15,6,1),$ which corresponds, after appropriate transfer of tens, to the number $1771561.$ This is $11^{6}.$ A first conclusion is therefore that these equalities are not mere coincidences and that there is a map taking $1+^nC_1+^nC_2+ \cdots+1=2^{n}=(1+1)^{n}$ to $(11_{10})^{n}.$ This map is expected to combine the coefficients, $(1,\;^nC_1,\;^nC_2,\; \cdots,\;1),$ of the Pascal triangle to form a whole number having the coefficients as the digits of the number (for $n<5$) or form the number after appropriate transfer of tens (for $n\geq5$). In order to define this map in its generality we shall first generalize the Pascal triangle.\\
\ \\

We consider $n,y \in \N \cup\{0\}$ and the coefficients $(1,\;^nC_1y,\;^nC_2y^2, \cdots,y^n)$ of the finite binomial expansion of $(x+y)^{n}.$ For different choices of $n,$ the corresponding triangle is $$1$$ $$1\;\;\;y$$ $$1\;\;\;2y\;\;\;y^{2}$$ $$1\;\;\;3y\;\;\;3y^{2}\;\;\;y^{3}$$ $$1\;\;\;4y\;\;\;6y^{2}\;\;\;4y^{3}\;\;\;y^{4}$$ $$1\;\;\;5y\;\;\;10y^{2}\;\;\;10y^{3}\;\;\;5y^{4}\;\;\;y^{5}$$ $$\cdots \cdots \cdots \cdots \cdots$$ $$\cdots \cdots \cdots \cdots \cdots$$ $$1\;\;\;^{n}C_{1}y\;\;\;^{n}C_{2}y^{2}\;\;\;^{n}C_{3}y^{3}\cdots^{n}C_{r}y^{r}\cdots y^{n}$$ $$\cdots \cdots \cdots \cdots \cdots$$ We shall refer to this as the \textit{Newtonian triangle} and denote it as $T(y).$ Its build-up formula may be seen as $y (^{n-1}C_{r-1}y^{r-1})+^{n-1}C_{r}y^{r}=^{n}C_{r}y^{r},\;\; \forall \; r \in \N,$ which becomes familiar when $y=1.$ In order to get a handle on our extension of the Pascal triangle we consider the Newtonian triangle for $y=2.$ In this case the $2$nd row is $(1,4,4),$ which corresponds to the number $144=12^{2},$ the $3$rd row is $(1,6,12,8)$ corresponding to the number $1728=12^{3},$ the $4$th row is $(1,8,24,32,16)$ corresponding to the number $20736=12^{4},\;etc.$ We shall therefore say that the number $20736$ \textit{digitally corresponds} to the row $(1,8,24,32,16),$ and vice-versa. We shall denote the \textit{digital-correspondence} map by $\delta:\N^{n+1}\rightarrow \N$ whose restriction to the subset $\{(1,\;^nC_1y,\;^nC_2y^2, \cdots,y^n):n,y\in \N\cup\{0\}\}$ of $\N^{n+1}$ is given as $$\delta(1,\;^nC_1y,\;^nC_2y^2, \cdots,y^n)=1\;\;^{n}C_{1}y\;\;^{n}C_{2}y^{2}\;\; \cdots\;\;y^{n},$$ where the right hand side is viewed as a whole number, whether tens are transferred (when $n\geq5$ or $y\neq1$) or not (when $n<5$ and $y=1$).\\
\ \\

The truth behind our observations that the whole number $\delta(1,\;^nC_1y,\;^nC_2y^2, \cdots,y^n)$ is always a power of $n$ may be formalized for any row $N(y,n):=(1,\;^nC_1y,\;^nC_2y^2, \cdots,y^n)$ in the Newtonian triangles. Here a natural number having an $n$th root in $\N,$ for some $n=2,3,4, \cdots,$ shall be called \textit{exact.}
\ \\
\ \\

{\bf 2.1. Lemma.}  \textit{Let $y,n \in \N\cup\{0\}$ and define $f_{n}(y)=\delta(N(y,n)).$ Each $f_{n}(y)\in \N$ and is exact of power $n.$ Every exact number in $\N$ is of the form $f_{n}(y).$}

{\bf Proof.}  We know that $(x+y)^{n}= x^{n}+(^{n}C_{1}y)x^{n-1}+(^{n}C_{2}y^{2})x^{n-2}+ \cdots + y^{n},$ so that, considering $x$ as the \textit{base of numeration} on both sides, we have $$(1y_{x})^n=(1\;\;^{n}C_{1}y\;\;^{n}C_{2}y^{2}\;\; \cdots\;\;y^{n})_{x}.$$ That is, $$(1y_{x})^n=\delta(N(y,n))_{x}\cdots \cdots \cdots(*)$$ as two equal numbers in base $x.\hspace{0.1in}\Box$
\ \\

The above Lemma shall be employed in Theorems $3.1,\;4.2$ and $4.3$ in the following form.

{\bf 2.2. Corollary.} \textit{Let $n \in \N$ be fixed and let $\mathfrak{E}_{n}$ be the collection of all exact rationals of power $n,$ explicitly given as $$\mathfrak{E}_{n}=\left\{\begin{array}{ll} \{\varepsilon^{n}:\varepsilon \in \mathbb{Q}^{+}\}, & \mbox{if}\;\;
n \in 2\mathbb{N},\\
\{\varepsilon^{n}:\varepsilon \in \mathbb{Q}\}, & \mbox{if}\;\;
n \in \mathbb{N} \setminus 2\mathbb{N}.
\end{array}\right.$$ Then the set $\mathfrak{E}_{n}$ is in a one-to-one correspondence with the set $\{f_{n}(y): y \in \mathbb{Q}\}.$}

{\bf Proof.} Define the map $\sigma $ as $\sigma(\varepsilon^{n})=f_{n}(y),$ where $\varepsilon \in \mathbb{Q},$ then $$\sigma: \mathfrak{E}_{n} \rightarrow \{f_{n}(y): y \in \mathbb{Q}\},$$ which, by Lemma $2.1,$ is a one-to-one correspondence. $\hspace{0.1in}\Box$
\\
\\

{\bf 2.3. Remarks on $f_{n}$}.\\
\ \\

$(1.)$  \textbf{On all exact rationals}: It may be seen, from the left side of $(*),$ that $f_{n}(y)=(10+y)^{n},$ as earlier envisaged in the case of $y=1.$ This polynomial form for $f_{n}$ allows us to extend its domain to all $y \in \mathbb{Q},$ giving \textit{only} \textit{all} exact rationals.
\ \\
\ \\

$(2.)$  \textbf{On Corollary 2.2}: A proof of Corollary $2.2$ which is independent of Lemma $2.1$ may also be given as follows. Define $\rho: \{f_{n}(y): y \in \mathbb{Q}\} \rightarrow \mathfrak{E}_{n}$ as $\rho(f_{n}(y)):=\varepsilon^{n},$ with $\varepsilon=10+y,\;y \in \mathbb{Q}.$ Clearly, $\rho$ is a one-to-one correspondence and $\rho=\sigma^{-1}.$
\ \\
\ \\

$(3.)$  \textbf{On general Diophantine equations}: Our focus is to discuss the contribution of $f_{n}(y)$ to Diophantine equations, which we may generally write as $$A_{1}\alpha_{1}^{n_{1}}+A_{2}\alpha_{2}^{n_{2}}+A_{3}\alpha^{n_{3}}_{3}+ \cdots + A_{p}\alpha^{n_{p}}_{p}=B\beta^{m},$$ for some constants $A_{i} \in \mathbb{Q},$ $n_{i} \in \N$ and unknowns $\alpha_{i} \in \mathbb{Q},\;i=1,2, \cdots,p.$ This translates, in our context, to studying $$A_{1}f_{n_{1}}(y_{1})+A_{2}f_{n_{2}}(y_{2})+A_{3}f_{n_{3}}(y_{3})+ \cdots + A_{p}f_{n_{p}}(y_{p})=Bf_{m}(y),$$ for some $y,y_{i} \in \mathbb{Q},\;i=1,2, \cdots,p.$ A particular example is when $A_{i}=1$ and $n_{i}=n$ with $p=2,$ which is the defining equation of $FLT.$ That is, $$f_{n}(y_{1})+f_{n}(y_{2})=f_{n}(y_{3}),$$ for $y_{1}\neq y_{2}\neq y_{3}.$ It is necessary to illustrate the depth of insight of this formulation of Diophantine equations by tackling a formidable problem.\ \\

We shall therefore illustrate our method with the problems of Pythagorean triples and $FLT.$ In our context, these two problems are simultaneously captured by studying the possible values of $y \in \mathbb{Q}$ for which $$Q_{n-1,a}(y):=f_{n}(y+a)-f_{n}(y),$$ $\;y \in \mathbb{Q},\;a \in \mathbb{Q} \setminus \{0\},\;n \in \N,$ is the digital-correspondence of some $N(y_{0},n),y_{0} \in \mathbb{Q}.$ We have set $y_{1}=y_{0},\;y_{2}=y\;\mbox{and}\;y_{3}=y+a,$ in $f_{n}(y_{1})+f_{n}(y_{2})=f_{n}(y_{3})$ above to arrive at the equation $Q_{n-1,a}(y):=f_{n}(y+a)-f_{n}(y)$.\\
\ \\
\ \\

Our approach is then to investigate, among other things, the \textit{reason} for the existence of (rational) Pythagorean triples (in \S $3.$), which we then employ to seek \emph{Fermat's triples,} if they exist (in \S $4.$).
\ \\
\ \\
\ \\
\ \\
\ \\
\ \\
\ \\

\indent {\bf \S 3. Pythagorean triples in the context of Newtonian triangles.} Lemma $2.1$ clearly says that $f_{2}(y),\forall\; y \in \N \cup \{0\},$ (indeed $\forall \; y \in \mathbb{Q}$) is a perfect-square in $\mathbb{Q}$ and that every perfect-square in $\mathbb{Q}$ is some $f_{2}(y).$ Hence the study of $f_{2}(y)$ translates to studying the digital-correspondence of the second rows, $N(y,2),$ of the Newtonian triangles, $T(y),$ for different values of $y.$ In this case $Q_{1,a}(y)=(2a)y+a(20+a).$ The following result may be seen as a purely algebraic and rational proof of the existence of Pythagorean triples and of the truth of Pythagoras' theorem for rationals. It establishes, in our context, that some of the values of $Q_{1,a}(y)$ appear in the list of the digital-correspondences of $N(y,2).$
\ \\
\ \\

{\bf 3.1. Theorem.} \textit{Let $a \in \mathbb{Q}\backslash \{0\}.$ Then there exist $y \in \mathbb{Q}$ for which $Q_{1,a}(y)$ is a perfect-square. That is, $Q_{1,a}(y)=\delta(N(y+b,2)),$ for some $y \in \mathbb{Q},\;b \in \mathbb{Q} \setminus \{0,a\}.$}

{\bf Proof.} Since $Q_{1,a}(y)$ is a linear polynomial in $y$ we substitute $y = \alpha_{2}x^2+\alpha_{1}x+\alpha_{0},\;x\in \mathbb{Q},$ where the
values of $\alpha_{2},\alpha_{1},\alpha_{0} \in \mathbb{Q}$ are yet to be known, into $Q_{1,a}(y)$ in order to consider $Q_{1,a}(y)$ for a candidate in the list of values of $\delta(N(y+b,2)).$ That is, $$Q_{1,a}(y)=Q_{1,a}(x)=(2a\alpha_{2})x^2+(2a\alpha_{1})x+a(2\alpha_{0}+20+a)$$ and, for it to be a complete square of a non-zero rational, we must have $Q_{1,a}(x) \equiv (px+q)^2$ for \textit{all} $p,q,x\in \mathbb{Q}.$ The choice of $y$ and the above identity are informed by the one-to-one correspondence in Corollary $2.2,$ with $n=2.$

This identity gives $\alpha_{2} = \frac{1}{2a}p^2,\alpha_{1} = \frac{1}{a}pq$ and $\alpha_{0} = \frac{q^{2}-a(20+a)}{2a},$ each of which belongs to $\mathbb{Q}$ uniquely, for \textit{every} $p,q\in \mathbb{Q}.$ Hence $$y = \left(\frac{p^2}{2a}\right)x^2+\left(\frac{pq}{a}\right)x +
\left[\frac{q^2-a(20+a)}{2a}\right]$$ is the required $y$ in $\mathbb{Q}.$ Indeed, the discriminant of $Q_{1,a}(x)$ vanishes exactly
when $\alpha_{2} = \frac{1}{2a}p^2,\alpha_{1} = \frac{1}{a}pq$ and $\alpha_{0} = \frac{q^{2}-a(20+a)}{2a}.$ \hspace{0.1in}$\Box$ \\
\ \\
\ \\

The conclusion of Theorem $3.1$ is that, for every $x,p,q\in \mathbb{Q}$ and $a \in \mathbb{Q} \setminus\{0\},$ the rational solutions, $y,$ to the equation $Q_{1,a}(y)=\delta(N(y+b,2))$ exist and are given as $y = \left(\frac{p^2}{2a}\right)x^2+\left(\frac{pq}{a}\right)x +
\left[\frac{q^2-a(20+a)}{2a}\right].$ The converse question to this result is that: \textit{if this $y$ is a given rational solution of $Q_{1,a}(y)=\delta(N(y+b,2)),$ does it imply that $x\in \mathbb{Q}?$} This question is addressed in the following Theorem.\\
\ \\
\ \\

{\bf 3.2. Theorem.} \textit{Let $p,q\;\mbox{and}\;a$ be as in the proof of Theorem $3.1,$ with $p \neq0.$ Every rational solution $y = \left(\frac{p^2}{2a}\right)x^2+\left(\frac{pq}{a}\right)x +
\left[\frac{q^2-a(20+a)}{2a}\right]$ of $Q_{1,a}(y)=\delta(N(y+b,2))$ corresponds to a rational value of $x.$}

{\bf Proof.} It is clear, from Theorem $3.1,$ that, if $x,p,q \in \mathbb{Q}$ and $a \in \mathbb{Q} \setminus \{0\},$ then the given $y$ is a solution of $Q_{1,a}(y)=\delta(N(y+b,2))$ and $y \in \mathbb{Q}.$ Conversely, let the given $y$ be a rational solution of $Q_{1,a}(y)=\delta(N(y+b,2))$ and let $(\alpha,\beta,\gamma)$ be a rational Pythagorean triple with $\alpha<\beta<\gamma.$ (That is, $(Q_{1,a}(y))^{\frac{1}{2}}<(f_{n}(y))^{\frac{1}{2}}<(f_{n}(y+a))^{\frac{1}{2}}$). Then $Q_{1,a}(y)=\alpha^{2}.$ This gives, $2ay+a(20+a)=\alpha^{2}.$ That is, $y=\frac{1}{2a}[\alpha^{2}-a(20+a)].$ Hence, $\frac{1}{2a}[\alpha^{2}-a(20+a)]=y = \left(\frac{p^2}{2a}\right)x^2+\left(\frac{pq}{a}\right)x +\left[\frac{q^2-a(20+a)}{2a}\right]$ which reduces to a quadratic equation in $x$ given as $p^{2}x^{2}+2pqx+(q^{2}-\alpha^{2})=0,$ with $p\neq0,$ which is necessary in order to find $x.$ The solution of this quadratic is $x=\frac{-q\pm \alpha}{p}\in \mathbb{Q}.$ \hspace{0.1in}$\Box$ \\
\ \\
\ \\

{\bf 3.3. Remarks on Theorem $3.1$}.\\
\ \\

$(1.)$  \textbf{On the coefficients of $y$}: Observe that it is necessary and sufficient for all $\alpha_{i},\;i=0,1,2,$ to be rational in order to always have $y \in \mathbb{Q}.$ The polynomial $Q_{1,a}(x)$ is \textit{always} a perfect-square of members of $\mathbb{Q}\setminus \{0\},$ whatever the value of $x$ in $\mathbb{Q}.$ A closer look at Theorem $3.1$ therefore reveals a very important conclusion that: \textit{in order to justify the identity used, between $Q_{1,a}(x)$ (which is always a perfect-square in $\mathbb{Q} \setminus \{0\}$) and $(px+q)^2,$ $p$ and $q$ must necessarily assume all values in $\mathbb{Q},$ and not just $\ 'some'$ values in $\mathbb{Q}.$} This observation, which is the core of the method of Theorem $3.1,$ shall be needed when considering rational Pythagorean triples and the non-zero rational solutions (if any) of $u^{n}+v^{n}=w^{n},$ for $n>2.$ See also $(1.)$ of Remark $(4.4).$\\
\ \\
\ \\

$(2.)$  \textbf{On the constant $b$}: Now that we have a general expression for $y \in \mathbb{Q}$ that explains the existence of Pythagorean triples, we may compute the constant $b \in \mathbb{Q} \setminus \{0,a\}$ in $Q_{1,a}(y)=\delta(N(y+b,2))$ as follows: $Q_{1,a}(y)=\delta(N(y+b,2)),$ for rational $y,$ $\Longleftrightarrow$ $y^{2}+(20+2b-2a)y+(b^{2}+10b+100-a^{2}-20a)=0$ has a perfect-square discriminant $\Longleftrightarrow$ the quadratic $2a^{2}-2ba+10b,$ in $a,$ has zero discriminant $\Longleftrightarrow$ $b=20.$\\
\ \\

A complete list of all \textit{rational} Pythagorean triples is therefore possible \textit{without} necessarily having to generate them from the basic example of the triple $(3,4,5).$ (See $[5.]$).
\ \\
\ \\

{\bf 3.4. Corollary.} \textit{Let $a,p,q,x \in \mathbb{Q}$ with $a\neq0.$ The general expression for any rational Pythagorean triple is then $(\alpha,\beta,\gamma)=$}
$$\left\{\begin{array}{ll} (\left(px+q\right),\left(\frac{p^2}{2a}\right)x^2+\left(\frac{pq}{a}\right)x +
\left(\frac{q^2-a^{2}}{2a}\right),\left(\frac{p^2}{2a}\right)x^2+\left(\frac{pq}{a}\right)x +
\left(\frac{q^2+a^{2}}{2a}\right)), & \mbox{if}\;
\alpha < \beta<\gamma,\\
(\left(\frac{p^2}{2a}\right)x^2+\left(\frac{pq}{a}\right)x +
\left(\frac{q^2-a^{2}}{2a}\right),\left(px+q\right),\left(\frac{p^2}{2a}\right)x^2+\left(\frac{pq}{a}\right)x +
\left(\frac{q^2+a^{2}}{2a}\right)), & \mbox{if}\;
\beta <\alpha<\gamma.
\end{array}\right.$$

{\bf Proof.} We already know, from Theorem $3.1,$ that every Pythagorean triple in $\mathbb{Q}$ is
$$(\alpha,\beta,\gamma)=\left\{\begin{array}{ll} (\sqrt{Q_{1,a}(y)},\sqrt{f_{2}(y)},\sqrt{f_{2}(y+a)}), & \mbox{if}\;
\alpha < \beta<\gamma,\\
(\sqrt{f_{2}(y)},\sqrt{Q_{1,a}(y)},\sqrt{f_{2}(y+a)}), & \mbox{if}\;
\beta <\alpha<\gamma,
\end{array}\right.$$ where $y$ is as found in the Theorem. Computing each of these triples with the said $y$ gives the result$.\;\Box$
\ \\
\ \\

\textit{We arrive at the classical Diophantus's solution to the problem of primitive solutions to $\alpha^{2}+\beta^{2}=\gamma^{2},$ if we set $x=0$ in Corollary $3.4$ and clear the fractions.}
\ \\
\ \\

{\bf 3.5. Corollary.} \textit{Let $\alpha \in \mathbb{Q}.$ Then there are $\beta,\gamma \in \mathbb{Q}$ such that $(\alpha,\beta,\gamma)$ is a rational Pythagorean triple. That is, every rational number is a first element of some rational Pythagorean triple.}

{\bf Proof.} Every $\alpha \in \mathbb{Q}$ may be written as $\alpha=px+q$ for a choice of $p,q,x \in \mathbb{Q}.$ The values of $\beta$ and $\gamma$ may then be computed from Corollary $3.4$ for any $a \in \mathbb{Q} \setminus \{0\}.\;\Box$\\
\ \\

The particular cases of \textit{non-trivial,} \textit{primitive} and  \textit{integral} Pythagorean triples may be deduced from these Corollaries, which may themselves be extended to include the study of Pythagorean $n-$tuples. See $[1.],\;$ p. $76.$ The generality inherent in the use of Newtonian triangles is evident from the ease with which general Pythagorean triples are handled. We now partition all integral Pythagorean triples into disjoint classes.\\
\ \\

Let $\mathbb{P}$ denote the set of all rational Pythagorean triples and denote the subset consisting of integral ones by $\mathbb{P}_{\mathbb{Z}}.$ Let $\mathbb{P}_{m}=\{(\alpha,\beta,\gamma) \in \mathbb{P}_{\Z}:\;gcd(\alpha,\beta,\gamma)=m\},$ where $m \in \Z.$ Clearly $\mathbb{P}_{\Z}= \bigcup_{m \in \Z} \mathbb{P}_{m}.$ It may not be clear whether or not this is a disjoint union. This may be addressed by using an appropriate equivalence relation.
\ \\
\ \\

{\bf 3.6. Theorem.} \textit{The equality $\mathbb{P}_{\Z}= \bigcup_{m \in \Z} \mathbb{P}_{m}$ is a disjoint union.}

{\bf Proof.} Define a relation $\sim$ on members of $\mathbb{P}_{\Z}$ as $(\alpha_{1},\beta_{1},\gamma_{1}) \sim (\alpha_{2},\beta_{2},\gamma_{2})$ \textit{iff} $gcd(\alpha_{1},\beta_{1},\gamma_{1})=gcd(\alpha_{2},\beta_{2},\gamma_{2}).$ It is immediate that $\sim$ is an equivalence relation on $\mathbb{P}_{\Z}.$ It is also clear that each $\mathbb{P}_{m}$ is a typical equivalence class in $\mathbb{P}_{\Z}/ \sim.\;\Box$\\
\ \\
\ \\

It therefore follows that the set $\{\mathbb{P}_{m}:m \in \Z\}$ is a partition of $\mathbb{P}_{\Z}.$
\ \\
\ \\

{\bf 3.7. Remarks on Pythagorean triples and  Newtonian triangles.}\\
\ \\

$(1.)$  \textbf{On parametrization of rational Pythagorean triples}: We may as well use $f_{2}(\lambda y)$ in the manner in which $f_{2}(y+a)$ has been considered. The first result here is that, for every $\lambda \in \mathbb{Q} \setminus \{0,1\},$ we always have that $$f_{2}(\lambda y)=\lambda^{2}f_{2}(y)+R_{1,\lambda}(y),$$ where $R_{1,\lambda}(y)=20 \lambda (1-\lambda)+100(1-\lambda^{2}).$ It can readily be shown that $R_{1,\lambda}(y)=\delta(N(\xi y,2)),\;\xi \in \mathbb{Q} \setminus \{0,\lambda\},$ \textit{iff} $$y = \left[\frac{p^2}{20\lambda(1-\lambda)}\right]x^2 +
\left[\frac{pq}{10\lambda(1-\lambda)}\right]x +
\left[\frac{q^2-100(1-\lambda^2)}{20\lambda(1-\lambda)}\right]$$ for
all $p,q,x \in \mathbb{Q}.$ This gives another outlook to Corollary $3.4.$
\ \\
\ \\

$(2.)$ \textbf{On equivalence classes of Pythagorean triples}: The function $$h: \mathbb{P}_{\Z}\rightarrow \Z$$ given as $h(\alpha,\beta,\gamma)=gcd(\alpha,\beta,\gamma),\;\forall\;(\alpha,\beta,\gamma) \in \mathbb{P}_{\Z},$ is well-defined and constant-valued on each $\mathbb{P}_{m}.$ It will be interesting to get the dependence of $m$ on the parameters of the triples in Corollary $3.4.$ That is, to derive a function $$\vartheta:\Z^{3}\times (\Z \setminus \{0\})\rightarrow \Z$$ given as $m=\vartheta(p,q,x,a)=gcd(\left(px+q\right),\left(\frac{p^2}{2a}\right)x^2+\left(\frac{pq}{a}\right)x +
\left(\frac{q^2-a^{2}}{2a}\right),\left(\frac{p^2}{2a}\right)x^2+\left(\frac{pq}{a}\right)x +
\left(\frac{q^2+a^{2}}{2a}\right)),$ where $(\left(px+q\right),\left(\frac{p^2}{2a}\right)x^2+\left(\frac{pq}{a}\right)x +
\left(\frac{q^2-a^{2}}{2a}\right),\left(\frac{p^2}{2a}\right)x^2+\left(\frac{pq}{a}\right)x +
\left(\frac{q^2+a^{2}}{2a}\right)) \in \mathbb{P}_{\Z},$ as this will put results on $h$ and $\mathbb{P}_{m}$ in proper perspectives. It may therefore be useful to note, from Corollary $3.4,$ that every $(\alpha,\beta,\gamma) \in \mathbb{P},$ with $\alpha<\beta<\gamma,$ (respectively, $\beta<\alpha<\gamma$) may be reduced to the (\textit{Diophantine}) form $(\alpha, \frac{\alpha^{2}-a^{2}}{2a},\frac{\alpha^{2}+a^{2}}{2a}),$ (respectively, $(\frac{\alpha^{2}-a^{2}}{2a},\alpha,\frac{\alpha^{2}+a^{2}}{2a})$), (where $\alpha:=px+q$ of Corollary $3.4$) for any $a \in \mathbb{Q}\setminus \{0\}.$ The well-known case of $\vartheta \equiv 1$ follows from here. The above \textit{Diophantine form} of the Pythagorean triples gives a compact expression for the result of Corollary $3.4$ and may be further discussed in the light of \textit{Hall's matrices,} $[5.].$
A step towards the derivation of an explicit expression for the function, $$\vartheta:\Z^{3}\times (\Z \setminus \{0\})\rightarrow \Z,$$ is to note, from the remark following Corollary $3.4,$ that $\vartheta(p,q,x,a)=1$ at $x=0.$ We may then write $\vartheta(p,q,x,a)=1+x\tau(p,q,x,a),$ where $\tau:\Z^{3}\times (\Z \setminus \{0\})\rightarrow \Z.$\\
\ \\
\ \\

$(3.)$ \textbf{On rings and modules of Newtonian triangles}: Let $n \in \N$ be fixed and consider the set $\mathfrak{N}(n):=\{N(y,n):\;y \in \Z\}.$ The operations $+$ and $\cdot,$ defined on members of $\mathfrak{N}(n),$ as $$N(y_{1},n)+N(y_{2},n):=N(y_{1}+y_{2},n)\;\;\mbox{and}\;\;N(y_{1},n) \cdot N(y_{2},n):=N(y_{1}y_{2},n),$$ respectively, convert $\mathfrak{N}(n)$ into a commutative ring with identity, $N(1,n),$ whose \textit{field of fractions} is $\{N(y,n):\;y \in \mathbb{Q}\}.$ The map $y\mapsto T(y)$ is a one-to-one correspondence between $\Z$ and $\mathfrak{N}(n),$ implying that $\mathfrak{N}(n)$ is indeed a Noetherian ring whose \textit{ideal} structure is exactly as in $\Z.$ If, in addition to these operations above, we define $\alpha N(y,n):=N(\alpha y,n)\;\alpha,y \in \Z,$ then $\mathfrak{N}(n)$ becomes a \textit{$\Z-$module.} These properties on $\mathfrak{N}(n)$ are inherited by the set $\mathfrak{T}_{\Z},$ of all Newtonian triangles, $T(y),\;y \in \Z,$ leading to the requirements that, for $y_{1},\;y_{2},\;y,\;\alpha \in \Z,$ $$T(y_{1})+T(y_{2}):=T(y_{1}+y_{2}),\;T(y_{1}) \cdot T(y_{2}):= T(y_{1}y_{2}),\;\mbox{and}\;\alpha T(y):=T(\alpha y).$$ In this formulation, \textit{the Pascal triangle, $T(1),$ is the (multiplicative) identity of the Noetherian ring $\mathfrak{T}_{\Z}$} while the \textit{functor,} $T,$ may be seen to be both \textit{covariant} and \textit{contravariant} on $\Z.$ The ring and module structures of $\mathfrak{T}_{\Z}$ are yet to be studied.
\ \\

In the light of our success on Pythagorean triples above, we are encouraged to consider the \textit{original} $FLT.$
\ \\
\ \\
\ \\
\ \\
\ \\
\ \\
\ \\

\indent {\bf \S 4. Fermat's Last Theorem in the context of Newtonian triangles.} The consideration of each $f_{n}(y),\;n>2,$ is essentially the study of the other rows, after the $2$nd, in each of the Newtonian triangles. Following in the direction of our method in \S 3., we compute the corresponding polynomial, $Q_{n-1,a}(y),\;n>2,$ which is then sought in the list of digital-correspondences to $N(y,n).$
\ \\
\ \\

{\bf 4.1. Lemma.} \textit{Let $a \in \mathbb{Q} \setminus \{0\}.$ Then  $Q_{n-1,a}(y)=nay^{n-1}+\frac{n(n-1)}{2!}(a^{2}+20a)y^{n-2}+\frac{n(n-1)(n-2)}{3!}(a^{3}+30a^{2}+300a)y^{n-3}+ \cdots +(a^{n}+10na^{n-1}+ \cdots +10^{n-1}na),$ for all $n \in \N,y \in \mathbb{Q}.$}

{\bf Proof.} Compute $f_{n}(y+a)-f_{n}(y).\hspace{0.1in}\Box$

In seeking a position for every $Q_{n-1,a}(y),\;n>2,$ in the list of digital-correspondence to $N(y,n)$ we make the following eye-opening observation on $Q_{2,a}(y).$
\ \\
\ \\

{\bf 4.2. Theorem.} ($cf.$ Euler's proof in $[3.],\;p.\;39.$) \textit{There does \textit{not} exist any $y \in \mathbb{Q}$ for which $Q_{2,a}(y)$ is a perfect cube. That is, $$Q_{2,a}(y)\neq\delta(N(y+b,3)),$$ $\forall\;y \in \mathbb{Q},\;b \in \mathbb{Q} \setminus \{0,a\}.$}

{\bf Proof.}  We assume the contrary and proceed as in Theorem $3.1.$ If the polynomial $Q_{2,a}(y) = 3ay^2+(3a^{2}+60a)y+(a^{3}+30a^{2}+300a)$ is to be a perfect-cube in $\mathbb{Q},$ there must exist $y=\alpha_{3}x^3+\alpha_{2}x^2+\alpha_{1}x+\alpha_{0}\in \mathbb{Q}$ with $\alpha_{3},\alpha_{2},\alpha_{1},\alpha_{0},x\in \mathbb{Q}$, such that, after substituting $y$ into $Q_{2,a}(y),$ the resulting polynomial, $Q_{2,a}(x),$ in $x$ and of degree six, would be identical to $(px^2+qx+r)^3,$ for all $p,q,r\in \mathbb{Q}.$ The choice of $y$ and the above identity are informed by the one-to-one correspondence in Corollary $2.2,$ with $n=3.$

By making this substitution and comparing the coefficients we arrive at seven relations, namely:
$\disp 3a\alpha^2_{3}=p^3,\;\;$ $2a\alpha_{2}\alpha_{3} = p^2q,\;\;$ $a(2\alpha_{1} \alpha_{3}+\alpha_{2}^2)=p^2r+pq^2,\;\;$\\
$6a(\alpha_{0} \alpha_{3}+\alpha_{1} \alpha_{2})+60a \alpha_{3}+3a^{2} \alpha_{3} = 6pqr+q^3$, $a(2\alpha_{0} \alpha_{2}+\alpha_{1}^2)+20a\alpha_{2}+ a^{2}\alpha_{2}=pr^2+q^2r$, $2\alpha_{0} \alpha_{1}+60a\alpha_{1}+3a^{2}\alpha_{1}=3qr^2$ and $\alpha^{2}_{0}+(3a^2+60a)\alpha_{0}+(a^{3}+30a^{2}+300a)=r^3,$ from which we are expected to find the rational constants $\alpha_{3},\alpha_{2},\alpha_{1}$ and $\alpha_{0}$ in terms of $p,\;q$ and $r$.  A consideration of the first
three and last relations give, if $p\neq0$ is assumed: $\disp
\alpha_{3}=\sqrt{\frac{1}{3a}p^3},\;\;\;$ $\disp \alpha_{2} =
\frac{p^{2}q\sqrt{3a}}{2a\sqrt{p^{3}}},\;\;\;$ $\disp
\alpha_{1}=\frac{(4p^2r+pq^2)}{8a}\left(\sqrt{\frac{3a}{p^{3}}}\right)\;\;\;$ and $\;\;\disp
\alpha_{0}=\frac{-3a^{2}-60a\pm\sqrt{9a^{4}+356a^{3}+3480a^{2}-1200a+4r^3}}{2}.$\\
These relations imply that $y \notin \mathbb{Q},$ if we use $(1.)$ of Remarks $(3.3).\hspace{0.1in}\Box$\\
\ \\

We may as well consider the use of $(px+q)^{6}$ instead of $(px^{2}+qx+r)^{3}$ in the proof of Theorem $4.2.$ However the use of $(px^{2}+qx+r)^{3}$ accommodates more generality than $(px+q)^{6},$ since not all quadratics are completely factorisable over $\mathbb{Q}.$ In any of these options the deduced expressions for $\alpha_{i},\;i=0,1,2,3$ do not satisfy the remaining three of the seven relations. A closer look at the proof reveals that this disorder in the identity, $Q_{2,a}(x)\equiv(px^{2}+qx+r)^{3},$ is primarily due to the disparity in the number of terms in $Q_{2,a}(x)$ (which is \textit{seven}) and the number of unknowns in the coefficients of $y$ (which is \textit{four}). There is no way to match these two numbers when $n>2,$ like what we have in the case of $n=2$ in Theorem $3.1,$ where there are three terms in $Q_{1,a}(x)$ and \textit{exactly} three unknowns in $y=\alpha_{2}x^{2}+\alpha_{1}x+\alpha_{0}.$\ \\

The method of proof of Theorem $4.2$ may be formalized in the following version of the \textit{original} $FLT.$
\ \\
\ \\

{\bf 4.3. Theorem.} \textit{Let $a \in \mathbb{Q} \setminus \{0\},\;n >2.$ Then there does \textit{not} exist any $y\in \mathbb{Q}$ for which $Q_{n-1,a}(y)$ is an exact rational of power $n.$ That is, $$Q_{n-1,a}(y)\neq \delta(N(y+b,n)),$$ $\forall\;y\in \mathbb{Q},\;b \in \mathbb{Q} \setminus \{0,a\}.$}

{\bf Proof.} We substitute $y=\alpha_{n}x^{n}+\alpha_{n-1}x^{n-1}+ \cdots +\alpha_{1}x+\alpha_{0}$ into $Q_{n-1,a}(y)$ in Lemma $4.1$ and observe that the only choices to be made of each $\alpha_{k},\;k=0,1,2, \cdots,n>2,$ for $Q_{n-1,a}(x)$ to be a digital-corresponding of some $N(y_{0},n),$ would involve extraction of roots, since powers of $y$ must have been computed in the process of substitution. This leads, via $(1.)$ of Remarks $(3.3),$ to the conclusion that $y \notin \mathbb{Q}.\hspace{0.1in}\Box$\\
\ \\

The above method of proof shows that a structural reason for the non-existence of Fermat's triples is because, in seeking a position for $Q_{n-1,a}(y)$ among the values of $\delta(N(y,n)),$ every substituted $y$ into $Q_{n-1,a}(y)$ must be raised to some powers, thereby introducing extraction of roots when coefficients of $y$ are later sought. The exception to this is in the cases of $n=1,2,$ where $Q_{n-1,a}(y)$ are the constant and linear polynomials, respectively. This explains the existence of rational triples, $(u,v,w),$ satisfying the Diophantine equations $u+v=w$ (when $n=1$ in $u^{n}+v^{n}=w^{n}$) and $u^{2}+v^{2}=w^{2}$ (when $n=2$ in $u^{n}+v^{n}=w^{n}$). A Galois equivalence of this reason has also been exploited in the next section. It is noted that no extra condition on $n,$ other than the original requirement of $n \in \Z$ and $n>2,$ was used to prove $FLT.$\\
\ \\
\ \\

{\bf 4.4. Remarks}.\\
\ \\

$(1.)$  \textbf{On the significance of the constant $a \in \mathbb{Q}\setminus \{0\}$}: Corollary $3.5$ reveals that every $a \in \mathbb{Q}\setminus \{0\}$ leads to a rational Pythagorean triple, while only \textit{some} $a \in \mathbb{Q}\setminus \{0\}$ gives the integral Pythagorean triples. The same may be deduced from the consideration of the non-zero rational and integral solutions of other Diophantine equations, say, $u^{3}-v^{3}=w^{2}.$ Indeed, substituting $y=\alpha_{1}x+\alpha_{0}$ into $Q_{2,a}(y),$ which, when identical with $(px+q)^{2},\forall\;p,q,x \in \mathbb{Q},$ gives $\alpha_{1}=\frac{p}{\sqrt{3a}},\;\alpha_{0}=\frac{2pq-(3a^{2}+60a)\alpha_{1}}{6a\alpha_{1}},$ we see that the non-zero rational solutions of $u^{3}-v^{3}=w^{2}$ exist only when $a=\frac{k^{2}}{3},\forall\;k \in \Z\setminus \{0\},p,q,x \in \mathbb{Q},\;p\neq0,$ while the non-zero integral solutions exist only when $a=\frac{1}{3},\;\mbox{for some}\;p,q,x \in \Z,\;p\neq0.$ It therefore follows that the non-zero rational constant $a,$ in $Q_{n-1,a}(y),$ measures the \textit{distribution} and \textit{density} of solutions of Diophantine equations, when they exist. This may be further explored\\
\ \\
\ \\

$(2.)$ \textbf{On unique factorization}: The method of this paper is to fix $n-$power of two arbitrary non-zero rationals, say $\alpha^{n}$ and $\beta^{n},$ and then seek for the possibility of a third one, $\gamma^{n},$ such that $\alpha^{n}+\beta^{n}=\gamma^{n},$ with $\alpha\beta\gamma\neq0.$ In this approach any two of $\alpha^{n},\beta^{n}$ and $\gamma^{n}$ may be fixed. However, our choice of $f_{n}(y+a)=f_{n}(y)+Q_{n-1,a}(y)$ over and above the other possibility of $f_{n}(y+a)+f_{n}(y)=P_{n,a}(y),$ which leads to the study of the  polynomial, $P_{n,a},$ of degree $n,$ is informed by the non-zero rational solutions of $\alpha^{2}+\beta^{2}=\gamma^{2}$ which, if considered in the light of $P_{n,a},$ will lead us outside the base field of $\mathbb{Q}.$ Indeed, considering any example of the \textit{Pythagorean triples,} say $(3,4,5),$ it is advisable, based on our approach, to use $5$ and $3$ to seek for $4$ by factorising the \textit{difference of two squares} $5^{2}-3^{2}$ as $5^{2}-3^{2}=(5-3)(5+3)=(2)(8)=(2)(2)(4)=4^{2}$ or to use $5$ and $4$ to seek for $3$ by factorising the \textit{difference of two squares} $5^{2}-4^{2}$ as $5^{2}-4^{2}=(5-4)(5+4)=(1)(9)=3^{2}$ than to use $3$ and $4$ to seek for $5$ by factorising the \textit{sum of two squares} $3^{2}+4^{2}$ as $3^{2}+4^{2}=3^{2}-(4i)^{2}=(3-4i)(3+4i)= \mid 3+4i \mid^{2}=5^{2},$ which, in the process, leads outside the base field of $\mathbb{Q}.$\\
\ \\

Thus, since factorisation in a fixed base field is the first step at extracting indices out of a number (and now, out of a polynomial), we have settled for the considerations of $Q_{n-1,a}(y)$ (which is the difference $f_{n}(y+a)-f_{n}(y)$), while we hope that the polynomials $P_{n,a}$ will be of immense use in aspects of  \textit{number theory} allowing the employment of the field $\mathbb{Q}(i)=\{a+bi \in \mathbb{C}:\;a,b \in \mathbb{Q}\}$ of \textit{gaussian numbers.} With the above approach we bypass the intricate manipulations involving \textit{unique factorization} in \textit{quadratic fields.} Other properties of the polynomials $f_{n}$ and $Q_{n-1,a},$ beyond their present use in the proof of $FLT,$ may also be studied.\\
\ \\
\ \\

$(3.)$ \textbf{On non-rational Pythagorean triples}: Our present approach in \S $3$ suggests the study of \textit{non-rational} Pythagorean triples in quadratic fields, $\mathbb{Q}(\sqrt[n]{\sigma}),$ (where $\sigma$ is an $n$th root-free rational number), in fields, $\mathbb{F}_{p},$ of prime characteristics and in fields, $\mathbb{Q}_{p},$ of \textit{p-adic} numbers. The significance of the constant $b=20$ in the present field of $\mathbb{Q},$ as derived in Theorem $3.1,$ or as may be derived in any other number field, is still unknown.\\
\ \\
\ \\
\ \\
\ \\
\ \\
\ \\
\ \\

\indent {\bf \S 5. Galois groups of Fermat Polynomials.} The \textit{original} Fermat's Last Theorem does \textit{not} translate to the investigation of solvability of the Galois group, $Gal(Q_{n-1,a}),$ of the polynomials $Q_{n-1,a},$ as it is always expected in the application of Galois theory to polynomials, but to the investigation of the values assumed by the order, $\mid Gal(Q_{n-1,a}) \mid,$ of $Gal(Q_{n-1,a}),$ as we shall show shortly. This approach around the \textit{Fermat polynomials,} $Q_{n-1,a},$ when combined with Theorem $4.3,$ gives the Galois group version of the original claim of Pierre de Fermat ([$3.$], p. $3$). The results of this section may also be used to deduce the nature of the roots of $Q_{n-1,a}(y)=0,\;\mbox{when}\;n>2.$\\
\ \\
\ \\

Let $L/K$ be a \textit{field extension.} We know that the \textit{degree}, $[L:K],$ of the extension satisfies $[L:K]=1$ \textit{iff} $L=K.$ If the extension is, in addition, \textit{normal} and \textit{separable} we conclude that the Galois group, $Gal(L/K),$ of the extension is the \textit{trivial} group. Now if $Gal(L/K)$ is the Galois group of a polynomial $f \in K[y],$ also written as $Gal(f)$ where $L$ is a splitting field of $f$ over $K,$ then $\mid Gal(L/K) \mid=1$ \textit{iff} $f$ has \textit{all} its roots in $K.$ That is, $\mid Gal(L/K)\mid=1$ \textit{iff} $f$ is completely reducible over $K.$ This observation may now be formalised.
\ \\
\ \\

\indent {\bf 5.1 Lemma.} \textit{$f \in K[y]$ is completely reducible over $K$ \textit{iff} $\mid Gal(f)\mid=1.$}

{\bf Proof.} Let $L$ be a splitting field of $f$ over $K.$ Then $L$ is a normal finite extension of $K$ and $\mid Gal(f)\mid=1$ \textit{iff} $[L:K]=1$ \textit{iff} $L=K.$ This means that $f$ is completely reducible over $K.\;\;\Box$
\ \\

Let $L$ be a splitting field of $f$ over $K.$ We shall call a polynomial $f \in K[y]$ \textit{incomplete with respect to $L/K$} whenever it has a linear factor in $L[y]$ which is not in $K[y].$ An \textit{opposite} to the above Lemma is therefore possible.\\
\ \\
\ \\

\indent {\bf 5.2 Lemma.} Let $L/K$ be a field extension. \textit{$f \in K[y]$ is incomplete with respect to $L/K$ \textit{iff} $\mid Gal(f)\mid\neq1.$}

{\bf Proof.} $(y-\alpha) \mid f(y)$ (for some $\alpha \in L\setminus K$) \textit{iff} $K[\alpha]$ is a splitting field of $f$ over $K$ \textit{iff} $[K[\alpha]:K]=2$ \textit{iff} $[N:K]\geq [K[\alpha]:K]=2\neq1$(where $N$ is the normal closure of $K[\alpha]$) \textit{iff} $\mid Gal(f)\mid=[N:K]\neq1.\;\Box$
\ \\

We may now study the \textit{Fermat polynomials,} $Q_{n-1,a},$ in the light of these Lemmas.
\ \\
\ \\

\indent {\bf 5.3 Theorem.} \textit{Each $Q_{n-1,a},$ with $n>2,\;a \in \mathbb{Q} \setminus \{0\},$ is an incomplete member of $\mathbb{Q}[y]$ with respect to any field extension of $\mathbb{Q}.$}

{\bf Proof.} We show that \textit{$\mid Gal(Q_{n-1,a}) \mid \neq 1,$ for all $n>2,\;a \in \mathbb{Q} \setminus \{0\}.$} Let $\alpha_{1},\alpha_{2}, \cdots ,\alpha_{n-1} \in \mathbb{C}$ be the roots of the monic polynomial $q_{n-1,a}:=\frac{1}{na}Q_{n-1,a},$ then, by the fundamental theorem of algebra, $$q_{n-1,a}(y)=(y-\alpha_{1})(y-\alpha_{2}) \cdots (y-\alpha_{n-1})=y^{n-1}-s_{1}y^{n-2}+s_{2}y^{n-3}+ \cdots +(-1)^{n-1}s_{n-1},$$ where $s_{1}=\alpha_{1}+\alpha_{2}+ \cdots + \alpha_{n-1}=\frac{-(n-1)}{2!}(a+20),$ $s_{2}=\alpha_{1}\alpha_{2}+\alpha_{1}\alpha_{3} \cdots + \alpha_{n-2}\alpha_{n-1}=\frac{(n-1)(n-2)}{3!}(a^{2}+30a+300),$ $\cdots ,s_{n-1}=\alpha_{1}\alpha_{2} \cdots \alpha_{n-1}=\frac{(-1)^{n-1}}{n}(a^{n-1}+10na^{n-2}+ \cdots +10^{n-1}n)$ are \textit{non-vanishing} elementary symmetric polynomials.
\ \\

Now let $L$ be a splitting field for $q_{n-1,a}$ over $\mathbb{Q}(s_{1},s_{2}, \cdots, s_{n-1}).$ Since the characteristics of $\mathbb{Q}$ is zero we conclude, from Theorem $10.10$ of $[7.],$ p. $178,$ that $Gal(q_{n-1,a})=S_{n-1}.$ Hence $Gal(Q_{n-1,a})=S_{n-1},$ because $\alpha_{1},\alpha_{2}, \cdots ,\alpha_{n-1}$ are also the roots of $Q_{n-1,a}.$ Thus $\mid Gal(Q_{n-1,a})\mid=(n-1)!.$ Since it is known that $(n-1)!=1$ \textit{iff} $n=1$ or $n=2,$ we therefore have that \textit{$\mid Gal(Q_{n-1,a}) \mid \neq 1,$ for all $n>2,\;a \in \mathbb{Q} \setminus \{0\}.$} $\Box$
\ \\

It is convenient to set $S_{n-1}$ as the trivial group, for $n=1,$ so that the popular choice of $0!$ as $1$ is justified. The complete treatment of the case of $n=2$ is contained in Theorem $2.1$ of $[9.].$
\ \\
\ \\

\indent {\bf 5.4 Remark.}\\

\indent $(1.)$ It follows therefore that, for $n>2,$ the normal closure of any splitting field of $Q_{n-1,a}$ over $\mathbb{Q}$ cannot be $\mathbb{Q}$ itself.  This is in contrast to the situation for $n=1,2.$
\ \\
\ \\

\indent $(2.)$ Also, since the proof of Theorem $5.3$ computes the group $Gal(Q_{n-1,a}),$ for \textit{all} $n \in \N,$ as $S_{n-1},$ whose order is $(n-1)!,$ we may therefore conclude that \textit{an underlying reason the equation $x^{n}=y^{n}+z^{n}$ has solutions in non-zero rationals only when $n=1$ (which follows from the field structure of $\mathbb{Q}$) and $n=2$ (as established in Theorem $3.1$), is because only $0!$ ($=\mid Gal(Q_{0,a})\mid,$ when $n=1$ in $\mid Gal(Q_{n-1,a})\mid$) and $1!$ ($=\mid Gal(Q_{1,a})\mid,$ when $n=2$ in $\mid Gal(Q_{n-1,a})\mid$) give the value $1$ among all $(n-1)!,$ $n \in \N.$} See the paragraph before Remarks $4.4$ for an equivalence of this reason.
\ \\
\ \\

\indent $(3.)$ \textbf{On Wiles-Taylor's proof of $FLT$}: It is expected that a profound theory would emerge out of the reconciliation of the modern theory of numbers, as has been put to use in $[13.],$ with the properties of the polynomials, $Q_{n-1,a},$ of the present paper. Indeed it would be interesting to link the analysis of Newtonian triangles to the \textit{Shimura-Taniyama-Weil} conjecture and the results of \textit{Diophantine geometry}.\\
\ \\
\ \\
\ \\
\ \\
\ \\
\ \\
\ \\

\indent {\bf \S 6. Arithmetic groups of Diophantine curves.} It is clear from above that $Q_{1,a}(x)=(2a)x+a(20+a)$ and that the non-zero rational points, $(x,y),$ on the \textit{Pythagorean} curve $$P_{a}:\;y^{2}=Q_{1,a}(x)=(2a)x+a(20+a)$$ are given as $(x,y)=(p^{2}\frac{z^{2}}{2a}+pq \frac{z}{a}+[\frac{q^{2}-a(20+a)}{2a}],pz+q),$ with $p,q,a \in \mathbb{Q} \setminus \{0\},z \in \mathbb{Q}.$ Define the non-empty set $G(P_{a})\subset \mathbb{Q}\times\mathbb{Q}$ as $G(P_{a}):=$ $$\{(x,y) \in \mathbb{Q}^{2}:x=p^{2}\frac{z^{2}}{2a}+pq \frac{z}{a}+[\frac{q^{2}-a(20+a)}{2a}],\;y=pz+q,\forall\;p,q,a \in \mathbb{Q} \setminus \{0\},z \in \mathbb{Q} \}$$ on which we define a binary operation as follows:
\ \\

Set $(x_{1},y_{1}),\;(x_{2},y_{2}) \in G(P_{a})$ as $(x_{1},y_{1})=(p^{2}_{1}\frac{z^{2}}{2a}+p_{1}q_{1} \frac{z}{a}+[\frac{q^{2}_{1}-a(20+a)}{2a}],p_{1}z+q_{1})$ and $(x_{2},y_{2})=(p^{2}_{2}\frac{z^{2}}{2a}+p_{2}q_{2} \frac{z}{a}+[\frac{q^{2}_{2}-a(20+a)}{2a}],p_{2}z+q_{2}),$ where $p_{i},q_{i},a \in \mathbb{Q} \setminus \{0\},\;z \in \mathbb{Q},\;i=1,2.$ We set $(x_{1},y_{1})\cdot(x_{2},y_{2}):=(x,y),$ where $$x=(p_{1}p_{2})^{2}\frac{z^{2}}{2a}+(p_{1}p_{2})(q_{1}q_{2}) \frac{z}{a}+[\frac{(q_{1}q_{2})^{2}-a(20+a)}{2a}]$$ and $$y=(p_{1}p_{2})z+(q_{1}q_{2}).$$ The following result then becomes immediate.\\
\ \\
\ \\

\indent {\bf 6.1 Theorem.} \textit{$(G(P_{a}),\cdot)$ is an abelian group whose identity element is given as $\textbf{1}=(\frac{z^{2}}{2a}+\frac{z}{a}+[\frac{1-a(20+a)}{2a}],z+1),$ with the inverse, $(x,y)^{-1},$ of every element, $(x,y)\in G(P_{a}),$ as $(x,y)^{-1}=((p^{-1})^{2}\frac{z^{2}}{2a}+(pq)^{-1} \frac{z}{a}+[\frac{(q^{-1})^{2}-a(20+a)}{2a}],p^{-1}z+q^{-1}).$}\\

\indent {\bf Proof.} We verify the well-known axioms of an abelian group. $\Box$
\ \\

It is known that, for each $a \in \mathbb{Q} \setminus \{0\},$ if $\mathbb{P}^{1}(\mathbb{Q})$ is the set of rational points of the projective $1-$space, then $G(P_{a})\simeq\mathbb{P}^{1}(\mathbb{Q})$ ($cf.\;\;[6.],$ Theorem $A.4.3.1$), so that each of the groups, $G(P_{a}),$ may be seen as a concrete realization of $\mathbb{P}^{1}(\mathbb{Q}).$\\
\ \\
\ \\

The method outlined above for the Pythagorean curve may be employed to compute the arithmetic group of any Diophantine curve. According to Theorem $4.3,$ the set $G(F_{a}):= \{(x,y)\in \mathbb{Q}^{2}:y^{n}=Q_{n-1}(x),\;n>2\},\;a \in \mathbb{Q}\setminus \{0\},$ consisting of \textit{non-zero} rational solutions of the \textit{Fermat} curve, $y^{n}=Q_{n-1}(x),\;n>2,$ is empty. For another example, the Diophantine curve attached to the non-zero rational solutions of $\alpha^{3}-\beta^{3}=\gamma^{2}$ is $y^{2}=Q_{2,a}(x).$ That is, the curve is $$E_{a}:\;y^{2}=(3a)x^{2}+(3a^{2}+60a)x+(a^{3}+30a^{2}+300a),$$ $a \in \mathbb{Q} \setminus \{0\}.$ The non-zero rational points on $E_{a}$ are then $$(x,y)=(\frac{p}{\sqrt{3a}}z+\frac{2pq\sqrt{3a}-(3a^{2}+60a)p}{6ap},pz+q),$$ where $p,q,z \in \mathbb{Q},\;p,q\neq0$ and $a=\frac{k^{2}}{3},\forall k \in \Z \setminus \{0\},$ from which the group operation may now be defined. The group $G(E_{a})$ is infinite, and since the genus of the curve $E_{a}$ is $0,$ it is also another concrete realization of $\mathbb{P}^{1}(\mathbb{Q}).$ However, finite subgroup of $G(E_{a})$ may be constructed from restrictions on its members. See $[1.],\;p.\;255,$ for an example of this restriction.
\ \\

This approach may be seen to have the capability of treating all the finiteness theorems of Diophantine geometry by explicitly computing the arithmetic group of any Diophantine curve. See $[6.],\;p.\;viii$ for a list of these theorems.
\ \\
\ \\

\indent {\bf 6.2 Remark: On attitude to a proof of $FLT.$}\\

It is somewhat sad that no one expects any longer that an elementary proof of the Fermat's Last Theorem will ever emerge. This is the conclusion of Michael Rosen ($[11.]$), some few years after the long, indirect and very difficult proof of Andrew Wiles and Richard Taylor was given in $[12.]$ and $[13.].$ This is borne out of the fact that many mathematicians were glad that the simple-looking statement of the Theorem could at least be \textit{said to have been finally proved} in $1994,$ after about $358$ years of sustained attacks by the most brilliant of each generation. What is really more grieving is the fact that the \textit{Wiles-Taylor proof} buried the totality of both the Theorem and the expectations of the rich theory that has been anticipated to emerge from its eventual proof, thus lending credence to the thought that $FLT$ is an isolated result of the theory of numbers.
\ \\
\ \\

This is exactly what is meant when Rosen said: \textit{To the degree that they (i.e., the partial results which appeared over the course of the centuries and which attempted to shed light on $FLT$) deal strictly with $FLT$ and not with any broader class of problems, it is an unfortunate fact that they are now obsolete.} Our approach in this paper therefore brings out the missed opportunities of the last three centuries that would have led straight to an easy understanding of the entire landscape of Diophantine Analysis of Equations, had it not been overlooked repeatedly. Indeed, if the $FLT$ is the non-existence result of rational solutions of $u^{n}+v^{n}=w^{n},\;n=3,4,5, \cdots,$ the polynomials, $Q_{n-1,a},\;\mbox{and}\;P_{n,a},$ deduced from it in $\S4.$ (and others that may be deduced from other Diophantine equations) are worthy of an independent study, as done in \S$5.,$ and of potential application to a wide range of subjects, as shown in the present section.\\
\ \\

Our present approach has the added advantage in that it \textit{does not deal strictly with $FLT,$} but, as may be seen in the last two sections, it is applicable to a wide range of subjects in algebraic number theory.\\
\ \\
\ \\
\ \\
\ \\
\ \\
\ \\
\ \\

\indent {\bf \S 7. Direct consequences of the Fermat's Last Theorem.} Contrary to what some experts in the modern theory of numbers would want us to believe, that the truth of the Fermat's Last Theorem (FLT) has no single application (even within number theory!) ($[4.]$ and $[8.]$), we consider some direct consequences of the Theorem in the form of open problems in the fields of topology, number theory, ring theory and Galois theory, all of which are deduced from the outlook of the proof of the Theorem given above. Hints on how these problems could be resolved are also included. It is our modest conclusion that the absence of these problems in the aftermath of the $1994$ Wiles-Taylor's proof of FLT is due mainly to the absence of other viable approaches to FLT and not that the truth of FLT has no single consequence.\ \\
\ \\
\ \\
\ \\

\indent {\bf \S A. On non-rational Fermat triples.} It is well-known that there are several \textit{quadratic} fields between the fields $\mathbb{Q}$ and $R,$ or between $\mathbb{Q}$ and $\C.$ One way of generating these subfields of $\R$ or of $\C$ is by the computation of the splitting fields of polynomials in, say, $\mathbb{Q}[X]$ or $\mathbb{Q}[X_{1},\cdots,X_{m}].$ The following problems are proposed:\ \\
\ \\

$(a.)$ Which of these splitting fields over $\mathbb{Q}$ will uphold the truth of the $FLT?$ That is, on which subfields, $\mathbb{F},$ of $\R$ or $\C$ is $Q_{n-1,a}(y)\neq \alpha^{n}$ for any $y,\alpha \in \mathbb{F}.$ The cases of $\mathbb{F}=\mathbb{Q}(\sqrt[n]{\sigma}),\;\mathbb{F}_{p},\;\mathbb{Q}_{p}$ have earlier been posited in Remark $4.4.3.$\ \\
\ \\

$(b.)$ Which of the splitting fields of the fermat polynomials, $Q_{n-1,a}$ (as may be deuced from Theorem $5.3,$ would admit the truth of $FLT$ and why?\\\
\ \\

$(c.)$ What is the numerical value and significance of the constant $b \in \mathbb{F}$ in the equation $Q_{n-1,a}(y)=\delta(N(y+b,n)),$ in those fields $\mathbb{F}$ that \textit{do not} admit the truth of the $FLT.$ For the case of $n=2$ and $\mathbb{F}=\mathbb{Q}$ we already know, from Remark $3.3.2,$ that $b:=b_{2,\mathbb{Q}}=20.$\ \\
\ \\
\ \\
\ \\

\indent {\bf \S B. Correct generalization of $\alpha^{2}+\beta^{2}=\gamma^{2}.$} The question has always been asked whether $FLT$ was the right question to the generalization of the Babylonian results on the sum of two (rational) squares being written as a rational square. It has been posited $([2.])$ that the correct analogue to the generalization of $\alpha^{2}+\beta^{2}=\gamma^{2}$ to cubes is not to consider $\alpha^{3}+\beta^{3}=\gamma^{3},$ but to seek non-zero rational solutions to $\alpha^{3}+\beta^{3}+\gamma^{3}=\delta^{3},$ while the situation for fourth powers is $\alpha^{4}+\beta^{4}+\gamma^{4}+\delta^{4}=\zeta^{4}, \cdots.$\\
\ \\
\ \\

In short, the conclusion of K. Choi $[2.]$ is that if rational solutions of $$x^{n}_{1}+ \cdots +x^{n}_{k}=z^{n}$$ are sought, it is necessary to first have that $k\geq n,$ though \textit{no} specific way of attacking this observation was suggested by him or by Davis Wilson (See Diophantine Equations on the website of $WolframMathWorld$) other than to state some conjectures and list the following suggestive examples: $3^{2}+4^{2}=5^{2}$ (where $k=2=n$), $3^{2}+4^{2}+12^{2}=13^{2}$ (where $k=3>2=n$), $3^{2}+4^{2}+12^{2}+84^{2}=85^{2}$ (where $k=4>2=n$), $3^{3}+4^{3}+5^{3}=6^{3}$ (where $k=3=n$), $4^{4}+6^{4}+8^{4}+9^{4}+14^{4}=15^{4}$ (where $k=5>4=n$), $4^{5}+5^{5}+6^{5}+7^{5}+9^{5}+11^{5}=12^{5}$ (where $k=6>5=n$), $\cdots.$ It is clear that there may be other examples that would escape the above scheme. We believe that the prospect of the case $k\geq n$ above should not preclude the investigation of the existence, or otherwise, of rational solutions of $x^{n}_{1}+ \cdots +x^{n}_{k}=z^{n}$ for $k<n,$ though it may require more than $100$ pages if we are to expect a proof of the  Wiles-Taylor's magnitude (which was the case $k=2<n$) to address \textit{each(!)} of the cases $2\neq k<n$ and the new cases of $k\geq n.$\ \\
\ \\

\textit{ We now propose an approach to this study (of both $k\geq n$ and $k<n$) based on an observation already contained in the proofs of Theorems $4.2$ and $4.3.$}\ \\
\ \\

With $k=2<n=3,$ we already have the non-existence of rational solutions of $x^{3}_{1}+x^{3}_{2}=z^{3}$ as Theorem $4.2$ above. A second look at the proof of this Theorem (as explained in the paragraph following it) shows that the conclusion of the Theorem stems from ~``the disparity in the number of terms in $Q_{n-1,a}(x)$ (which is seven) and the number of unknowns in the coefficients of $y$ (which is four)." It was also reported that there was \textit{no way} to match these two numbers in the case $k=2<n=3,$ unless we increase the number of cubes being added. That is, unless we increase $k$ beyond $2.$ Indeed if $k=3=n,$ then $x^{3}_{1}+x^{3}_{2}+x^{3}_{3}=z^{3}$ may be recast as $z^{3}-x^{3}_{1}-x^{3}_{2}=x^{3}_{3},$ which translate (in the context of Newtonian triangles) to studying a cubic polynomial $R_{3,a,b},$ given as $$R_{3,a,b}(y):=f_{3}(y+a)-f_{3}(y+b)-f_{3}(y),$$ for $y,a,b \in \mathbb{Q},\;a\neq0,\;b\neq0.$ We then seek $y \in \mathbb{Q}$ for which $R_{3,a,b}(y)=\delta(N(y+c,3)),$ $c \in \mathbb{Q} \setminus \{0,a,b\},$ where the \textit{lacuna} noted in the proof of Theorem $4.2$ would have been filled due to the introduction of the \textit{new} term, $f_{3}(y+b).$\ \\
\ \\

It is clear, form this paper, how the above outlined approach for $k=3=n$ may be achieved for all $k=3\geq n$ and indeed for any $k\geq n,$ whenever $n$ is fixed in $\N.$ We need only refer to Lemma $4.1$ for orientation on the general situation of $k=2< n,$ which may itself be extended to the most general case of $k<n.$ \ \\
\ \\

\textit{The present problem, as outlined above, is a strong argument in favour of our methods of handling the FLT and in the complete understanding of the study of Diophantine equations.}\ \\
\ \\
\ \\
\ \\

\indent {\bf \S C. On Fermat metric.} We consider here a direct consequence of $FLT$ and fix the positive intger $n\geq3.$ It is already shown that the polynomial, $Q_{n-1,a}(y),$ of Lemma $4.1$ is $\neq \delta(N(y+b,n))$ as long as $y \in \mathbb{Q} \setminus \{0,a\},$ but that $Q_{n-1,a}(y)=\delta(N(y+b,n)),$ whenever $y \in \R \setminus \{0,a\},$ for any choice of $b \in \R.$ A serious question along this line of thought is how the topologies on the two fields of $\mathbb{Q}$ and $\R$ contribute to the above conclusions about $Q_{n-1,a}(y)$ and $\delta(N(y+b,n)),$ since we know that, in the \textit{Euclidean metric,} $\overline{\mathbb{Q}}=\R.$ However, there are other topologies on $\mathbb{Q}$ in whose metric the completion, $\overline{\mathbb{Q}},$ would not be $\R.$ We mention the well-known $p-$adic completion, $\overline{\mathbb{Q}}=\mathbb{Q}_{p}.$ It is still an open problem, included in \S $A.$ above, if $Q_{n-1,a}(y)=\delta(N(y+b,n)),$ for any $y \in \mathbb{Q}_{p} \setminus \{0,a\}.$ These and many other examples of topologies and metrics on the subsets, $\N,\;\Z,\;\mathbb{Q},\;\mathbb{Q}_{p}, \cdots,$ of $\R$ (or of $\C$) lead to the consideration of the following definition:\ \\
\ \\

{\bf 7(C).1 Definition.} \textit{Let $(\overline{X},\rho)$ be the completion of a metric space, $(X,\rho).$ The metric, $\rho,$ is called a fermat metric if whenever FLT holds in $(X,\rho)$ it also holds in $(\overline{X},\rho).$ We then refer to the pair $(X,\rho)$ as a fermat metric space.}\ \\
\ \\

In other words a \textit{fermat metric} is a metric $\rho:X \times X \rightarrow [0,\infty)$ for which $Q_{n-1,a}(y)\neq \alpha^{n},$ for all $y,\;\alpha \in (\overline{X},\rho).$ If $X=\N,\;\Z$ and $x,\;y \in X,$ we set $\rho$ as $\rho(x,y)=\mid x-y\mid,$ then$(X,\rho)$ is a fermat metric space while $(\mathbb{Q},\rho)$ is not.\ \\
\ \\

In the general situation of the above definition one would like to know if every metric on a fermat metric space is a fermat metric and which of the topologies on $X$ may be deduced from a fermat metric. Also, for which example of the set, $X,$ (whether finite, discrete, Baire, $\cdots$) is every metric a fermat metric? A description of the open sets, closed sets, accumulation of a set, interior of a set, base of the topology, $\cdots,$ in terms of the fermat polynomials, $Q_{n-1,a}(y),\;y \in X,$ will contribute richly to our understanding of \textit{polynomial-induced metrics.} An open problem in \S $A.$ is to know whether or not the $p-$adic metric is a fermat metric on $\mathbb{Q}.$\ \\
\ \\
\ \\
\ \\

\indent {\bf \S D. Galois theory of Fermat fields.} This section may be seen as an analytic continuation of the exploration in \S $1.$ above. Let $F_{i},\;1\leq i\leq r$ be a collection of subfields of $\R$ (or $\C$). We shall call any member of this collection a \textit{fermat field} whenever $FLT$ holds on it. \ \\
\ \\

{\bf 7(D).1 Definition.} \textit{Let $\mathbb{Q} \subseteq \mathbb{F}_{1}\subseteq \cdots \subseteq \mathbb{F}_{r}\subset \R$ (or $\C$) be an increasing collection of fields. We refer to the collection, $F_{i},\;1\leq i\leq m,\;m\leq r,$ as a collection of nested fermat fields of length $m$ whenever $(a.)$ $\mathbb{Q} \subseteq \mathbb{F}_{1}\subseteq \cdots \subseteq \mathbb{F}_{m}\subset \R$ (or $\C$) and $(b.)$ each $\mathbb{F}_{i},\;1\leq i\leq m,$ is a fermat field.}\ \\
\ \\

Some of the important questions on this definition are:\ \\
\ \\

$(a.)$ How many collection of nested fermat fields are there for each exponent $n\geq 3?$\ \\
\ \\

$(b.)$ Is there a relationship between the length of a nested fermat fields and each $n?$\ \\
\ \\

$(c.)$ In the general case of $\mathbb{Q} \subseteq \mathbb{F}_{1}\subseteq \cdots \subseteq \mathbb{F}_{r}\subset \R$ (or $\C$), at what field, $\mathbb{F}_{k},\;1\leq k\leq r,$ does $FLT$ holds for which it fails at $\mathbb{F}_{k+1}$ and what is the relationship of $k$ to $n?$\ \\
\ \\

$(d.)$ What are the properties of $\mathbb{F}_{k}$ and $\mathbb{F}_{k+1}$ in $(c.)$ above and how does the Galois groups, $Gal(\mathbb{F}_{t+1}/\mathbb{F}_{t}),\;t=1,2,3, \cdots,$ of the field extensions, $\mathbb{F}_{t+1}/\mathbb{F}_{t},$ contribute to these conclusions above?\ \\
\ \\

$(e.)$ Is $Gal(\mathbb{F}_{t+1}/\mathbb{F}_{t})$ in any relationship with $Gal(Q_{n-1,a})$ (which has been computed above to be $S_{n-1}$) or with $Gal(P_{n,a})$ (with $P_{n,a}$ as in Remark $4.4.2$)?\ \\
\ \\

$(f.)$ How does an arithmetic group (if non-empty) of any Diophantine curve contributes to all these open problems?\ \\
\ \\
\ \\
\ \\

We are hoping to attack some of these open problems in collaboration with others.\\
\ \\
\ \\
\ \\
\ \\
\ \\
\ \\

\indent {\bf References.}
\begin{description}
\item [{[1.]}] Andreescu, T., Andrica, D. and Cucurezeanu, I., \textit{An introduction to Diophantine equations.} Birkh$\ddot{a}$user-Verlag. $2010.$
\item [{[2.]}] Choi, K., A note on Fermat's Last Theorem- Was it a right question?\ \\
\textit{www.public.iastate.edu/$\sim$kchoi/fermat.htm}
\item [{[3.]}] Edwards, H. M., \textit{Fermat's Last Theorem.} Graduate Text in Mathematics. {\bf 50}. Springer-Verlag. $1984.$
\item [{[4.]}] Frey, G., The way to the proof of Fermat's Last Theorem,\ \\
\textit{www.backup.itsoc.org/review/05pl1.pdf,} $p.\;1-17,\;1997.$
\item [{[5.]}] Hall, A., Geneology of Pythagorean triples, \textit{Math. Gaz.} LIV $(1970),$ $377-9;$ XI ${\bf 15.}$
\item [{[6.]}] Hindry, M. and Silverman, J. H., \textit{Diophantine geometry, an introduction.} Graduate Text in Mathematics. Springer-Verlag, $2000.$
\item [{[7.]}] Howie, J. M., \textit{Fields and Galois Theory.} Springer Undergraduate Mathematics Series. Springer-Verlag. $2006.$
\item [{[8.]}] Mazur, B., Number theory as gadfly, \textit{American Mathematical Monthly,} August-September, $1991,$ $p.\;593-610.$
\item [{[9.]}] Oyadare, O. O., On the application of Newtonian triangles to decomposition theorems. Preprint.
\item [{[10.]}] Page, A., \textit{Algebra.} University of London Press Ltd. First Published $1947.$ Reprinted $1950.$
\item [{[11.]}] Rosen, M., Book Review of Paulo Ribenboim's \textit{Fermat's Last Theorem for Amateurs.} in \textit{Notices Amer. Math. Soc.} \textbf{47}. ($2000$). $474-476.$
\item [{[12.]}] Taylor, R. and Wiles, A. J., Ring theoretic properties of certain Hecke algebras. \textit{Ann. of Math.} {\bf 141}. $(1995).$ $553-572.$
\item [{[13.]}] Wiles, A. J., Modular elliptic curves and Fermat's Last Theorem. \textit{Ann. of Math.} {\bf 141}. $(1995).$ $443-551.$
\end{description}

\ \\
Department of Mathematics,\\
Obafemi Awolowo University,\\
Ile-Ife, $220005,$\\
NIGERIA.\\
URL: \textit{www.maths.oauife.edu.ng/profile/Oyadare}\\
\text{Email address: \textit{femi\_oya@yahoo.com}}\\
\text{Mobile: +234(0)7031816287}

\end{document}